\documentclass[a4paper,12pt,reqno]{amsart}

\usepackage{amsfonts}
\usepackage{amsmath}
\usepackage{amssymb}

\usepackage{mathrsfs}
\usepackage{hyperref}

\newtheorem{theorem}{Theorem}[section]

\newtheorem{proposition}[theorem]{Proposition}
\newtheorem{lemma}[theorem]{Lemma}

\newcommand\R{{\bf R}}

\newcommand\pa{\partial}
\newcommand\Rn{{\bf R}^n}

\makeindex             

\begin{document}

\title[Global well-posedness of the Kirchhoff equations and systems]
{Global well-posedness of the Kirchhoff equation and Kirchhoff systems}

\author[Tokio Matsuyama]{Tokio Matsuyama}
\address{
  Tokio Matsuyama:
  \endgraf
  Department of Mathematics
  \endgraf
  Chuo University
  \endgraf
  1-13-27, Kasuga, Bunkyo-ku, Tokyo 112-8551 
  \endgraf
  Japan
  \endgraf
  {\it E-mail address} {\rm tokio@math.chuo-u.ac.jp}
  }

\author[Michael Ruzhansky]{Michael Ruzhansky}
\address{
  Michael Ruzhansky:
  \endgraf
  Department of Mathematics
  \endgraf
  Imperial College London
  \endgraf
  180 Queen's Gate, London SW7 2AZ 
  \endgraf
  United Kingdom
  \endgraf
  {\it E-mail address} {\rm m.ruzhansky@imperial.ac.uk}
  }

\date{\today}

\maketitle

\begin{abstract}
This article is devoted to review the known results on  
global well-posedness for the Cauchy problem to the Kirchhoff equation
and Kirchhoff systems with small data. Similar results will be obtained for the 
initial-boundary value problems in exterior domains with compact 
boundary. Also, the known results on large data problems will be 
reviewed together with open problems. 
\end{abstract}


\section{Introduction}
\label{sec:1}
G. Kirchhoff proposed the equation 
\[
\pa^2_t u-\left(1+\int_\Omega |\nabla u(t,y)|^2\, dy \right)\Delta u=0 \quad 
(t\in {\bf R}, \, x\in \Omega)
\]
in his book of mathematical physics in 1883, as a model equation for transversal motion of 
the elastic string, where $\Omega$ is a domain in ${\bf R}^n$ 
(see \cite{Kirchhoff} and also \cite{Nishida}). 
Since then, it was the first that in 1940 Bernstein proved the existence of global in time 
analytic solutions on an interval of real line in his celebrated 
paper \cite{Bernstein}. 
As is easily seen, this equation has a Hamilton structure,  nevertheless 
it involves a mysterious 
problem {\em whether or not, one can prove the existence of time global solutions with 
large data in Gevrey class, $C^\infty$ class or standard Sobolev spaces.} 
The results on global well-posedness in the 
Sobolev spaces $H^{3/2}$, or $H^2$ with small data are well established in 
\cite{Callegari,Dancona1,Dancona2,Dancona3,Greenberg,Heiming,Kajitani,Ext-Matsuyama,
Rend-Matsuyama,Racke,Rzymowski,Yamazaki1,Yamazaki2}. There, the classes of small data consist of 
compactly supported functions, or more generally, they are characterised by some weight conditions 
(see \cite{Callegari,Dancona1,Dancona2,Dancona3,Greenberg}) or oscillatory integrals 
(see \cite{Kajitani,Ext-Matsuyama,Rend-Matsuyama,Rzymowski}). Recently, the present authors 
studied the global well-posedness 
for Kirchhoff systems with small data (see \cite{MR-Liouville}), and extended all the previous results. 
Here, the class of data in \cite{MR-Liouville} consists of Sobolev space $(H^1)^m$, $m$ being 
the order of system, and 
is characterised by some oscillatory integrals. 

The object in this article is to review the known results cited above, and open problems 
for the Kirchhoff equation or Kirchhoff systems. 
More precisely, let us consider the Cauchy problem for the 
Kirchhoff equation
of the form
\begin{eqnarray}\label{EQ:Kirchhoff}
& \partial^2_t u-\left(1+\displaystyle{\int_{{\bf R}^n}} 
|\nabla u(t,y)|^2\, dy \right)\Delta u=0, & \quad t\ne 0, \quad x\in {\bf R}^n,\\
& u(0,x)=f_0(x), \quad \pa_t u(0,x)=f_1(x), &\quad x \in {\bf  R}^n.
\label{EQ:Kirchhoff-initial}
\end{eqnarray}
Kirchhoff systems are of also interest. 
For $m$ column vector-valued function $U=U(t,x)$ on $\R\times \Rn$,
consider 
\begin{eqnarray}\label{EQ:S-Kirchhoff}
& D_t U=A(s(t),D_x) U, & \quad t\ne 0, \quad x\in {\bf R}^n,\\ 
& U(0,x)= {}^t (f_0(x),f_1(x),\cdots,f_{m-1}(x)), &\quad x \in {\bf  R}^n,
\label{EQ:S-initial}
\end{eqnarray}
where 
\[
D=i^{-1}\partial \quad (i=\sqrt{-1}),
\] 
$A(s,D_x)$ is an $m\times m$ pseudo-differential system
with a suitably smooth element for $s\in \R$ in a neighbourhood of 
$0$; $s(t)$ is a quadratic form 
\[
s(t)=\langle SU,U\rangle_{{\bf L}^2(\Rn)}
=\int_{\Rn} {}^t (SU(t,x))\overline{U(t,x)}\, dx, 
\] 
where $S$ is an $m\times m$ Hermitian matrix,
and we put 
$${\bf L}^2(\Rn)=(L^2(\Rn))^m.$$ 
Assume that system (\ref{EQ:S-Kirchhoff}) is 
strictly hyperbolic:
Characteristic polynomial of  
$D_t-A(s,D_x)$ has real and distinct roots 
$\varphi_1(s,\xi),\ldots,\varphi_m(s,\xi)$ for any $s$ in the domain of definition 
of the $A(s,\xi)$ and for any $\xi\in \Rn\backslash0$, i.e., 
\begin{equation}
\mathrm{det}(\tau I-A(s,\xi))=(\tau-\varphi_1(s,\xi)) \cdots 
(\tau-\varphi_m(s,\xi)).
\label{Kirchhoff strict hyperbolicity1}
\end{equation} 
Hence, for $\delta>0$ sufficiently small, we may assume that 
\begin{equation}\label{Kirchhoff strict hyperbolicity2}
\inf_{|\xi|=1,\, s \in[0,\delta]} 
|\varphi_j(s,\xi)-\varphi_k(s,\xi)|>0
\end{equation}
for $j \ne k$.
Needless to say, Kirchhoff systems (\ref{EQ:S-Kirchhoff}) cover 
the Kirchhoff equation (\ref{EQ:Kirchhoff}).\\

The approach of this paper yields new results already in the scalar
case of the classical Kirchhoff equation but, in fact, we are able
to make advances for coupled equations as well or, more generally, 
for Kirchhoff systems (see \S \ref{sec:2}).
To this end, Kirchhoff systems are of interest but present several 
major complications compared to the scalar case. 
First of all, even for the linearised system, 
it is much more difficult to find a suitable representation of solutions which
would, on one hand, work with the low regularity ($C^1$) of coefficients
while, on the other hand, allow one to obtain sufficiently good estimates for
solutions. Moreover, in the case of systems of higher order, it is impossible to
find its characteristics explicitly, and the geometry of the system or rather of
the level sets of the characteristics enters the picture. \\

We will present the three kinds 
of well-posedness considered in the literature, namely:
\begin{enumerate}
\item Kirchhoff equation with large data 

\item Kirchhoff equation with small data 

\item Kirchhoff systems with small data 

\end{enumerate}

The main idea behind items 2--3 is to approach the problem by 
developing the ``asymptotic integration'' method for the linearised equation to be able to control its solutions 
to the extent of being able to prove a priori estimates necessary for the handling of the fully nonlinear 
problem. \\

We recall the definition of Sobolev spaces. We denote by 
$H^\sigma(\Rn)=(1-\Delta)^{-\sigma/2}L^2(\Rn)$ for 
$\sigma\in \R$ the standard Sobolev spaces,
and their homogeneous version is 
$\dot{H}^\sigma(\Rn)=(-\Delta)^{-\sigma/2}L^2(\Rn)$.
We put $H^\infty(\Rn)=\displaystyle{\bigcap_{j=0}^\infty} H^j(\Rn)$.


\section{Model of Kirchhoff systems}
\label{sec:2}
In this section, for the convenience of the considerations, 
we shall give some examples of Kirchhoff systems 
(\ref{EQ:S-Kirchhoff}). These examples are applicable 
for our result. The first example is:\\

\noindent 
{\bf Example 2.1.} Let us consider the Cauchy problem 
to the Kirchhoff equations of higher order:
\begin{eqnarray*} 
& L\left(D_t,D_x,s(t)\right)u
\equiv D^m_t u+
\displaystyle{\sum_{j \le m-1,\vert \nu \vert+j=m}} 
b_{\nu,j}\left(s(t)\right) 
D^{\nu}_x D^{j}_t u=0, \\
& D^{k}_t u(0,x)=f_k(x), \quad 
k=0,1,\cdots,m-1. 
\end{eqnarray*}
Here the quadratic form $s(t)$ is given by 
\[
s(t)=\int_{\Rn}\sum_{|\beta|=|\gamma|=m-1}s_{\beta \gamma}D^\beta u(t,x)
\overline{D^\gamma u(t,x)}\, dx,
\]
where $\beta=(\beta_t,\beta_x)$, $\gamma=(\gamma_t,\gamma_x)$, 
$D^\beta=D_t^{\beta_t}D_x^{\beta_x}$ and 
$s_{\beta \gamma}=\overline{s_{\gamma\beta}}$. 
We assume that the symbol $L(\tau,\xi,s)$ of the differential operator 
$L(D_t,D_x,s)$ has real 
and distinct roots 
$\varphi_1(s,\xi),\ldots,\varphi_m(s,\xi)$ 
for $\xi \ne0$ and $0\le s \le \delta$ with $\delta>0$, i.e., 
\[
L(\tau,\xi,s)=(\tau-\varphi_1(s,\xi)) \cdots 
(\tau-\varphi_m(s,\xi)), 
\]
\[
\inf_{s \in[0,\delta],\, |\xi|=1} 
|\varphi_j(s,\xi)-\varphi_k(s,\xi)|>0
\]
for $j \ne k$.
By taking the Fourier transform in the space variables and introducing the vector 
\[
V(t,\xi)={}^T(\vert \xi \vert^{m-1} \widehat{u}(t,\xi),\vert \xi \vert^{m-2} D_t \widehat{u}(t,\xi),
\cdots,D^{m-1}_t \widehat{u}(t,\xi)),
\]
we reduce the problem to the  system 
\begin{eqnarray*}
D_t V
=& \left(
\begin{array}{cccc}
0 & 1 & \ldots & 0 \\ 
0 & 0 & \ddots & 0 \\
\vdots & \ddots & \ddots & 1 \\
-H_m(s(t),\xi) 
& -H_{m-1}(s(t),\xi) 
& \ldots 
& -H_1(s(t),\xi) 
\end{array}\right)|\xi|V
=A(s(t),\xi)V,
\end{eqnarray*}
where we put
\[
H_j(s(t),\xi)=\sum_{\vert \nu \vert=j} b_{\nu,m-j}(s(t)) (\xi/|\xi|)^\nu,
\quad (j=1,\ldots,m).
\]

As a new example of (\ref{EQ:S-Kirchhoff}), 
we can treat the completely coupled Kirchhoff equations of
the following type. \\

\noindent 
{\bf Example 2.2.} 
Let us consider the Cauchy problem to the coupled Kirchhoff equations 
of the following type:
\begin{eqnarray} \label{EQ:Example}
& \partial_t^2u
-a_1\left(1+\|\nabla u(t)\|_{L^2}^2+\|\nabla v(t)\|_{L^2}^2\right)
\Delta u+P_1(t,D_x)v=0, \\ 
& \partial_t^2v
-a_2\left(1+\|\nabla u(t)\|_{L^2}^2+\|\nabla v(t)\|_{L^2}^2\right)
\Delta v+P_2(t,D_x)u=0, \nonumber\\
& \partial_t^j u(0,x)=u_j(x), 
\quad \partial_t^j v(0,x)=v_j(x), \quad 
j=0,1,
\nonumber
\end{eqnarray}
for some second order homogeneous polynomials $P_1(t,D_x),P_2(t,D_x)$, and 
for some constants $a_1,a_2>0$ with $a_1\ne a_2$.
The quadratic form is given here by
\[
s(t)=\|\nabla u(t)\|_{L^2}^2+\|\nabla v(t)\|_{L^2}^2.
\] 
We assume that
\begin{eqnarray}\label{EQ:Ass-AP1}
& |\xi|^{-2}P_k(t,\xi)\in  \mathrm{Lip}_{\mathrm{loc}}(\R;L^\infty(\Rn\backslash0)), \\
& |\xi|^{-2}\partial_t P_k(t,\xi)\in  L^1(\R;L^\infty(\Rn\backslash0)),
\end{eqnarray}
for $k=1,2$, and that
\begin{equation}
\inf_{t\in \R,\, |\xi|=1} 
\left\{(a_1-a_2)^2 +4P_1(t,\xi)P_2(t,\xi)\right\}>0, 
\end{equation}
\begin{equation}\label{EQ:Ass-AP2}
\inf_{t\in \R,\, |\xi|=1} 
\left\{a^2_1 a^2_2-P_1(t,\xi)P_2(t,\xi)\right\}>0.
\end{equation}
By taking the Fourier transform in the space variables and introducing the vector 
\[
V(t,\xi)
={}^t(|\xi|\widehat{u}(t,\xi),\widehat{u}^\prime(t,\xi),
|\xi|\widehat{v}(t,\xi),\widehat{v}^\prime(t,\xi)),
\]
we rewrite (\ref{EQ:Example}) as a system 
\[
D_t V= \left(
\begin{array}{cccc}
0&-i|\xi|&0&0\\
ic_1(t)^2|\xi|&0& iP_1(t,\xi)|\xi|^{-1}&0\\
0&0&0& -i|\xi|\\
iP_2(t,\xi)|\xi|^{-1} &0& ic_2(t)^2|\xi|&0
\end{array}\right)V
=: A(s(t),\xi)V, 
\]
where 
\[
c_k(t)=\sqrt{a_k(1+s(t))}, \quad k=1,2.
\]
The four characteristic roots of the equation
$$\mathrm{det}(\tau I-A(s(t),\xi))=0$$ 
in $\tau$
are given by
\[
\varphi_{1,2,3,4}(s(t),\xi)
=
\pm \frac{|\xi|}{\sqrt{2}}
\sqrt{c_1(t)^2+c_2(t)^2\pm
\sqrt{\{c_1(t)^2-c_2(t)^2\}^2+4P_1(t,\xi) P_2(t,\xi)|\xi|^{-4}}}.
\]
Then it follows from (\ref{EQ:Ass-AP1})--(\ref{EQ:Ass-AP2}) that 
\[
\inf_{s \in[0,\delta],\, |\xi|=1} 
|\varphi_j(s,\xi)-\varphi_k(s,\xi)|>0
\]
for $j \ne k$.


\section{Known results and open problems with large data}
\label{sec:3}
In this section we shall review some known results on 
global well-posedness of the Kirchhoff equation with large data. 
Let us introduce some classes $\mathcal{N}_{\varphi}$ 
which ensure the global well-posedness for 
(\ref{EQ:Kirchhoff})--(\ref{EQ:Kirchhoff-initial}):
\[
\mathcal{N}_{\varphi}=\bigcup_{\eta>0} \left\{(f_0,f_1)\in H^1\times L^2:
\int_{\Rn} \left\{|\xi|^2|\hat{f}_0(\xi)|^2+|\hat{f}_1(\xi)|^2\right\} 
e^{\eta \varphi(|\xi|)}\, d\xi<\infty \right\}
\]
for some function $\varphi:[0,\infty)\to [0,\infty)$.
Here $\hat{f}(\xi)$ stands for the Fourier transform of $f(x)$. 

When $\varphi_0(|\xi|)=|\xi|$, $\mathcal{N}_{\varphi_0}$ coincides with 
analytic class $\mathcal{A}_{L^2}$ of $L^2$ type.
Bernstein proved that $\mathcal{A}_{L^2}$ well-posedness for $n=1$ (see 
\cite{Bernstein}). After him,
Arosio and Spagnolo discussed analytic solutions in 
higher spatial dimensions (see \cite{Arosio}). 
Also, we refer to the result of Kajitani and Yamaguti who proved 
analytic well-posedness for the degenerate Kirchhoff equation
(see \cite{Kajitani-Pisa}).
There are also quasi-analytic classes which ensure the global 
well-posedness. Here 
quasi-analytic classes are intermediate ones between the 
analytic class and $C^\infty$-class.
Nishihara found a weight $\varphi_1(|\xi|)=|\xi|/\log (e+|\xi|)$, and 
treated more general weight functions $\varphi_1$ with convexity 
condition (see \cite{Nishihara}). 
For the sake of simplicity, we treat only 
$\varphi_1(|\xi|)=|\xi|/\log (e+|\xi|)$. Then he 
proved that if the data belong to $\mathcal{N}_{\varphi_1}$,  
$H^\infty$-solutions exist globally. 
Ghisi and Gobbino generalised Nishihara class in the sense that 
the convexity condition is removed (see \cite{Ghisi}). 
Summarising the above known results, we have the inclusions:
\[
\mathcal{A}_{L^2} \subset \mathcal{N}_{\varphi_1} 
\subset \mathcal{G}^s
\subset H^\infty\times H^\infty \quad (s>1),
\]
where inclusions are all strict.
Here $\mathcal{G}^s$ is the Gevrey class of $L^2$ type and order $s$;  
\[
\mathcal{G}^s=\mathcal{N}_{|\xi|^{1/s}}, \quad (s\ge1).
\]
When $s=1$, $\mathcal{G}^1$ coincides with $\mathcal{A}_{L^2}$.
It should be noted that $\mathcal{G}^s$ well-posedness is still open.\\

We have other special classes which ensure the global well-posedness.
The first class is of Poho$\check{\rm z}$haev. He introduced the following class
(see \cite{Pohozhaev}):
\[
\mathcal{P}=\left\{(f_0,f_1)\in H^\infty\times H^\infty:
\limsup_{j\to\infty} 2^j \left(\|\Delta^{2^j+(1/2)}f_0\|^2_{L^2}
+\|\Delta^{2^j}f_1\|^2_{L^2}
\right)^{-\frac{1}{2^{j+2}}}>0
\right\}.
\]
The second one is of Manfrin, the class $\mathcal{B}_\Delta$
(see \cite{Manfrin-JDE}, and also \cite{Hirosawa}): 
we say that $(f_0,f_1)\in \mathcal{B}_\Delta$ if $(f_0,f_1)\in H^2\times H^1$ and 
there exists a positive sequence $\{\rho_j\}_{j\in {\bf N}}$, $\rho_j\to\infty$,
and a constant $\eta>0$ such that
\[
\sup_{j\in {\bf N}}\int_{|\xi|\ge \rho_j} \left\{
|\xi|^4|\hat{f}_0(\xi)|^2+|\xi|^2|\hat{f}_1(\xi)|^2
\right\} \frac{e^{\eta \rho^2_j/|\xi|}}{\rho^2_j}\, d\xi<\infty.
\] 
The inclusion among these classes are:
\[
\mathcal{A}_{L^2}\subset \mathcal{P} \subset \mathcal{B}_\Delta \subset 
H^2\times H^1,
\]
\[
\mathcal{P} \not\subset \mathcal{N}_{\varphi_1}, \quad \mathcal{G}^s
\not\subset \mathcal{B}_\Delta \quad (s>1).
\]
In conclusion, $\mathcal{G}^s$-well-posedness, $C^\infty$-well-posedness 
and $H^s$-well-posedness are still open. 
These problems are quite fascinating.


\section{Kirchhoff equation with small data}
\label{sec:4}
In this section we will review the recent results on global 
well-posedness for initial-boundary value problem to 
the Kirchhoff equation with small data. 
For an open set $\Omega$ in $\Rn$ with a smooth boundary
$\pa\Omega$, 
we consider the initial-boundary value problem:
\begin{eqnarray}\label{EQ:Equation}
& \pa^2_t u-
\left(1+\displaystyle{\int}_\Omega |\nabla u|^2\, dx\right)
\Delta u=0, &\quad t\in \R, \quad x\in\Omega,\\
& u(0,x)=f_0(x), \quad \pa_t u(0,x)=f_1(x), & \quad x\in \Omega,\\
& u(t,x)=0, & \quad t\in \R, \quad x\in \pa\Omega.
\label{EQ:Dirichlet}
\end{eqnarray}
Then we have:
\begin{theorem}[Matsuyama (\cite{Rend-Matsuyama})]
\label{thm:Base}
Let $n\ge1$. If $(f_0,f_1)\in Y(\Omega)$ and
\[
\|\nabla f_0\|^2_{L^2(\Omega)}+\|f_1\|^2_{L^2(\Omega)}
+|(f_0,f_1)|_{Y(\Omega)}\ll 1,
\]
then the initial-boundary value problem
(\ref{EQ:Equation})--(\ref{EQ:Dirichlet}) admits a
unique solution $u\in C(\R;H^{3/2}(\Omega)\cap H_0^1(\Omega))
\cap C^1(\R;H^{1/2}(\Omega))$, where 
\[
Y(\Omega):=\left\{(f,g)\in (H^{3/2}(\Omega)\cap H_0^1(\Omega))
\times H^{1/2}(\Omega)\, : \,
|(f,g)|_{Y(\Omega)}<+\infty \right\},
\]
with
\begin{eqnarray*}
 |(f,g)|_{Y(\Omega)} &=& \int_{-\infty}^{+\infty} 
\left\{
\left|\left(e^{i\tau H}H^{3/2}f,H^{3/2}f \right)_{L^2(\Omega)} \right|
+\left|\left(e^{i\tau H}H^{3/2}f,H^{1/2}g\right)_{L^2(\Omega)}\right| 
\right. \\
& & \left. 
+\left|\left(e^{i\tau H}H^{1/2}g,H^{1/2}g \right)_{L^2(\Omega)}
\right| 
\right\} 
\, d\tau.
\end{eqnarray*}
Here $(f,g)_{L^2(\Omega)}$ denotes the $L^2(\Omega)$-inner product
of $f$ and $g$.
\end{theorem}

When $\Omega=\Rn$, a similar space is considered by Kajitani  
and Rzymowski (see \cite{Kajitani,Rzymowski}). As to the prevous 
known results on the global well-posdness for small data, 
we refer to the results of Greenberg and Hu, D'Ancona and Spagnolo, and 
of Yamazaki 
(see \cite{Greenberg,Dancona1,Dancona2,Dancona3,Yamazaki}). 
Greenberg and Hu proved that 
if the data belong to $H^2(\R)\times H^1(\R)$ and 
have compact supports, $H^2(\R)$-solution exists globally. 
After them, D'Ancona and Spagnolo
obtained the similar results in higher spatial dimensions if the data 
belong to the weighted Sobolev space $H^2_\kappa (\Rn) \times 
H^1_\kappa (\Rn)$ for $\kappa\in(1,n+1]$, where we define
\[
H^s_\kappa(\Rn)=\{
f\in \mathcal{S}^\prime(\Rn):
\langle x \rangle^\kappa f \in H^s(\Rn)\}, \quad s\in\R.
\]
More general result was discussed by Yamazaki. She introduced 
a class $\mathcal{Y}_\kappa(\Rn)$ which ensures the global well-posedness.
More precisely, $(f_0,f_1)\in \mathcal{Y}_\kappa(\Rn)$ is said to be
$(f_0,f_1)\in H^{3/2}(\Rn)\times H^{1/2}(\Rn)$ such that
\[
\sum_{j,k=0}^1 \sup_{\tau\in \R}\langle \tau\rangle^\kappa
\left|\int_{\Rn} e^{i\tau|\xi|} \widehat{f}_j(\xi)\overline{\widehat{f}_k(\xi)}
|\xi|^{3-j-k}\, d\xi \right|<\infty.
\]
We have the inclusion among these classes as follows:
\[
H^2_\kappa (\Rn) \times H^1_\kappa (\Rn)
\subset \mathcal{Y}_\kappa (\Rn)
\subset Y(\Rn).
\]  
We notice that the first inclusion holds provided 
$\kappa\in (1,n+1]$, and 
the second one is valid for any $\kappa>1$.
In conclusion, our calss $Y(\Rn)$
is the most general which ensures global well-posedness 
for the (scalar) Kirchhoff equation. \\

As to the exterior problems, we have 
\begin{theorem}[Matsuyama (\cite{Ext-Matsuyama})] \label{prop-J}
Let $\Omega$ be a domain of $\Rn$ such that
$\Rn\setminus\Omega$ is compact and its boundary $\partial \Omega$
is $C^\infty$.
For $\sigma\ge0$ and $\kappa\in \R$,
let $H^\sigma_{\kappa,0}(\Omega)$ be the completion of
$C^\infty_0(\Omega)$
in the norm $\|\cdot\|_{H^\sigma_\kappa(\Omega)}$.
If $n\ge3$ and $\Rn\setminus \Omega$ is star-shaped
with respect to the origin,
then the inclusion
\[
H^{s_0+1}_{s(k),0}(\Omega)\times H^{s_0}_{s(k),0}(\Omega)
\subset Y_k(\Omega) \subset Y(\Omega)
\]
holds for any $s_0>(n+1)/2$, $s(k)>\max(n+1/2,k+n/2)$ and $k>1$.
\end{theorem}

The proof of Theorem \ref{prop-J} 
is based on the generalized Fourier transformation method. 
The cruicial tool in our argument is the asymptotic expansion of resolvent 
of $-\Delta$ around the origin in the complex plane. \\

Based on Theorems \ref{thm:Base}--\ref{prop-J}, we have: 
\begin{theorem}[Matsuyama (\cite{Ext-Matsuyama})] 
\label{thm:1-Main Theorem}
Let $\Omega,n,s_0,s(k)$ be as in Theorem \ref{prop-J}.
If
\[
f_0(x)\in H^{s_0+1}_{s(k),0}(\Omega), \quad
f_1(x)\in H^{s_0}_{s(k),0}(\Omega),
\]
and
\[
\|f_0\|_{H^{s_0+1}_{s(k)}(\Omega)}
+\|f_1\|_{H^{s_0}_{s(k)}(\Omega)}
\ll1,
\]
then the initial-boundary value problem
(\ref{EQ:Equation}) admits a
unique solution 
\[
u\in \bigcap_{j=0}^2 C^j(\R;H^{s_0+1-j}(\Omega)).
\]
\end{theorem}

Let us now make only a few short remarks to compare the result 
in Theorem \ref{thm:1-Main Theorem} with what is known. 
In the results of Heiming and Racke the data are imposed to be small in the 
weighted Sobolev spaces, and the supports of generalised Fourier 
transform of data are away from the origin (see \cite{Heiming,Racke}). 
However, our result covers 
the low frequencies of the generalised Fourier transform of data. 
Therefore, the statement of Theorem \ref{thm:1-Main Theorem} 
goes beyond \cite{Heiming,Racke}.

We should refer to the results of Yamazaki 
(see \cite{Yamazaki1,Yamazaki2}), who gave some sufficient conditions 
without any weight condition on data. That is, she assumed that the data 
belong to $W^{s,q}(\Omega)\times W^{s-1,q}(\Omega)$ for some 
$s>2$ and $q\in (1,2)$ depending on $n(\ge3)$, where $\Omega$ 
is a non-trapping domain in $\Rn$. 
We have an advantage of considering the classes of Theorem 
\ref{thm:1-Main Theorem}; it is more useful in the scattering problem 
rather than the ones in \cite{Heiming,Racke,Yamazaki1,Yamazaki2}. \\

As a final remark in this section, we have a new result on the Cauchy 
problem. Indeed, we can remove the additional class 
$\mathcal{Y}$ of data. This result will appear elsewhere. 


\section{Kirchhoff systems with small data}
\label{sec:5}
Recall Kirchhoff system (\ref{EQ:S-Kirchhoff}):
\begin{eqnarray*}
& D_t U=A(s(t),D_x) U, & \quad t\ne 0, \quad x\in {\bf R}^n,\\ 
& U(0,x)= {}^t (f_0(x),f_1(x),\cdots,f_{m-1}(x)), &\quad x \in {\bf  R}^n.
\end{eqnarray*}
Denote 
$${\bf H}^\sigma(\Rn)=(H^\sigma(\Rn))^m
$$ 
for $\sigma\in \R$.\\

Then we have:
\begin{theorem}[Matsuyama and Ruzhansky (\cite{MR-Liouville})]
\label{thm:S-Kirchhoff}
Let $n\ge1$. Suppose that $A(s,\xi)=(a_{jk}(s,\xi))_{j,k=1}^m$ is an $m\times m$ matrix 
positively homogeneous order one in $\xi$, whose entries $a_{jk}(s,\xi/|\xi|)$ are in 
$\mathrm{Lip}([0,\delta];L^\infty(\Rn\backslash0))$ for some $0<\delta \ll1$, 
and satisfies the strictly hyperbolic condition
(\ref{Kirchhoff strict hyperbolicity1})--(\ref{Kirchhoff strict hyperbolicity2}). 
If $U_0(x)\in {\bf L}^2(\Rn)\cap \mathcal{Y}(\Rn)$ satisfy 
\begin{equation}\label{EQ:small}
\|U_0\|^2_{{\bf L}^2(\Rn)}+|U_0|_{\mathcal{Y}(\Rn)}
\ll 1,
\end{equation}
then system (\ref{EQ:S-Kirchhoff})--(\ref{EQ:S-initial}) 
has a unique solution 
$U(t,x)\in C(\R;{\bf L}^2(\Rn)).$
In addition to (\ref{EQ:small}), if $U_0\in {\bf H}^1(\Rn)$,
then the solution $U(t,x)$ exists uniquely in the class 
$C(\R;{\bf H}^1(\Rn))\cap C^1(\R;{\bf L}^2(\Rn))$.
\end{theorem}

Precise definition of $\mathcal{Y}(\Rn)$ is as follows: 
$U_0={}^t(f_0,f_1,\ldots,f_{m-1})\in \mathcal{Y}(\Rn)$ if and only if  
$U_0\in (\mathcal{S}^\prime(\Rn))^m$ satisfies 
\begin{eqnarray*}
|U_0|_{\mathcal{Y}(\Rn)}:&=& \sum_{j,k=0}^{m-1}
\int^\infty_{-\infty}\left(
\int_{{\bf S}^{n-1}}\left|\int_0^\infty 
e^{i\tau \rho}
\widehat{f}_j(\rho\omega)\overline{\widehat{f}_k(\rho\omega)} 
\rho^n \, d\rho \right|\, d\sigma(\omega)\right)
\, d\tau\\
&<& +\infty,
\end{eqnarray*}
where $d\sigma(\omega)$ is the $(n-1)$-dimensional 
Hausdorff measure.\\

Let us compare Theorem \ref{thm:S-Kirchhoff} with what is known. 
Callegari and Manfrin introduced the following class (see \cite{Callegari}, 
and also \cite{Manfrin1}): 
\[
\mathcal{M}(\Rn)=\left\{U_0(x)={}^t (f_0(x),f_1(x),\ldots,f_{m-1}(x)):
|U_0|_{\mathcal{M}(\Rn)}<\infty \right\},
\]
where 
\[
|U_0|_{\mathcal{M}(\Rn)}=\sum_{k=0}^2\sum_{j=0}^{m-1}
\sup_{\omega\in {\bf S}^{n-1}} \int^\infty_0
\left|\partial^k_\rho \widehat{f}_j(\rho\omega)\right|^2 
\left(1+\rho^{\max\{n,2\}} \right)\, d\rho.
\] 
The inclusion among this class and ours is:
\[
(C_0^\infty(\Rn))^m\subset 
{\bf L}^1_2(\Rn)\cap {\bf H}^1_2(\Rn)\subset 
\mathcal{M}(\Rn) \subset 
\mathcal{Y}(\Rn),
\]
where denoting $L^1_2(\Rn)=\{f\in 
\mathcal{S}^\prime(\Rn):\langle x \rangle^2f\in L^1(\Rn)\}$ and 
$H^1_2(\Rn)=\{f\in
H^1(\Rn):\langle x \rangle^2 f\in L^2(\Rn)\}$, we put 
$${\bf L}^1_2(\Rn)
=(L^1_2(\Rn))^m,
\quad {\bf H}^1_\kappa(\Rn)
=(H^1_\kappa(\Rn))^m.
$$
Therefore, we can understand that Theorem \ref{thm:S-Kirchhoff} 
is the most general result in the framework of small data 
problem.\\

As another application of Theorem \ref{thm:S-Kirchhoff}, 
we consider the Cauchy problem for the second order equation of the form 
\begin{equation}\label{EQ:Spagnolo}
\partial^2_t u-\left(1+\int_{\Rn} |u(t,y)|^2 \, dy \right)\Delta u=0, 
\quad t\ne0, \quad x \in \Rn,
\end{equation}
with data 
\begin{equation} \label{EQ:Spagnolo-data}
u(0,x)=f_0(x), \quad \partial_t u(0,x)=f_1(x).
\end{equation}
In this particular case, the nonlocal term $s(t)$ is defined by 
\[
s(t)=\|u(t)\|^2_{L^2(\Rn)}.
\]
Introducing another class of data 
\[
\widetilde{\mathcal{Y}}(\Rn)
=\left\{(f_0,f_1)\in \mathcal{S}^\prime(\Rn)\times \mathcal{S}^\prime(\Rn):
|(f_0,f_1)|_{\widetilde{\mathcal{Y}}(\Rn)}<\infty\right\},
\]
where we put 
\begin{eqnarray*}
|(f_0,f_1)|_{\widetilde{\mathcal{Y}}(\Rn)}&=&\\
&& \sum_{j,k=0}^1
\int^\infty_{-\infty}\left(
\int_{{\bf S}^{n-1}}\left|\int_0^\infty
e^{i\tau \rho}
\widehat{f}_j(\rho\omega)\overline{\widehat{f}_k(\rho\omega)} 
\rho^{n-j-k} \, d\rho \right|\, d\sigma(\omega)\right)
\, d\tau,
\end{eqnarray*}
we have{\rm :}

\begin{theorem}[Matsuyama and Ruzhansky (\cite{MR-Liouville})] 
\label{thm:NONLOCAL}
Let $n\ge1$. For any 
$(f_0,f_1)\in (H^1(\Rn)\times L^2(\Rn))\cap \widetilde{\mathcal{Y}}(\Rn)$, 
the Cauchy problem (\ref{EQ:Spagnolo})--(\ref{EQ:Spagnolo-data}) has a unique solution 
$u\in \displaystyle{\bigcap_{k=0,1}}C^k(\R;H^{1-k}(\Rn))$, provided that 
\[
\|f_0\|^2_{L^2(\Rn)}+\|f_1\|^2_{\dot{H}^{-1}(\Rn)}+
|(f_0,f_1)|_{\widetilde{\mathcal{Y}}(\Rn)}\ll 1.
\]
\end{theorem}

When $n\geq 3$,
a similar result was obtained in Callegari and Manfrin 
and D'Ancona and Spagnolo (see \cite{Callegari,Dancona2}).  
However, 
the regularity of data in Theorem \ref{thm:NONLOCAL} is lower than 
that in  the previous results. 
It should be noted that Theorem \ref{thm:NONLOCAL} 
also covers low dimensions $n=1,2$, the case that remained open since 
Callegari and Manfrin, and D'Ancona and Spagnolo. \\


\section{Outline of proof of Theorem \ref{thm:S-Kirchhoff}}
The strategy of the proof of Theorem \ref{thm:S-Kirchhoff}
is to employ the Schauder-Tychonoff fixed point theorem 
via asymptotic integrations method.
Consider the {\it linear} Cauchy problem: 
\begin{eqnarray}\label{EQ:Eq}
&  D_t U=A(t,D_x)U, \quad (t,x)\in \R\times \Rn,\\ 
&  U(0,x)=U_0(x)={}^t \left(f_0(x),f_1(x),\ldots,f_{m-1}(x)\right)
\label{EQ:initial}
\end{eqnarray}
where $A(t,D_x)$ is a first order $m\times m$ pseudo-differential system,
with symbol $A(t,\xi)$ of the form
$A(t,\xi)=(a_{jk}(t,\xi))_{j,k=1}^m$. 
We assume that 
$\mathrm{det} (\tau I-A(t,\xi))=0$ has real and distinct roots 
$\varphi_1(t,\xi),\ldots,\varphi_m(t,\xi)$,
\begin{equation}\label{hyp1}
a_{jk}(t,\xi/|\xi|) \in \mathrm{Lip}_{\mathrm{loc}}(\R;L^\infty(\Rn\backslash0))
\quad \mathrm{and} \quad 
\partial_t a_{jk}(t,\xi/|\xi|) \in L^1(\R;L^\infty(\Rn\backslash0)).
\end{equation}

We prepare the next lemma. 
\begin{lemma}[c.f. Proposition 6.4 from \cite{Mizohata}] 
\label{Diagonalisation}
Let $A(t,\xi)$ be a symbol of differential operator $A(t,D_x)$
satisfying (\ref{hyp1}).
Then there exists a matrix $\mathcal{N}=\mathcal{N}(t,\xi)$ of homogeneous 
order $0$ in $\xi$ satisfying the following properties{\rm :} 

\noindent 
{\rm (i)} $\mathcal{N}(t,\xi) A(t,\xi/|\xi|)
=\mathcal{D}(t,\xi)\mathcal{N}(t,\xi)$, 
where 
\[
\mathcal{D}(t,\xi)
=\mathrm{diag} \left( \varphi_1(t,\xi/|\xi|),\ldots,
\varphi_m(t,\xi/|\xi|)\right);
\]

\noindent 
{\rm (ii)} $\displaystyle
{\inf_{\xi \in \Rn\backslash 0, t \in \R}}
\vert {\rm det} \, \mathcal{N}(t,\xi)) \vert>0;$

\noindent 
{\rm (iii)} $\mathcal{N}(t,\xi)\in \mathrm{Lip}_{\mathrm{loc}}
(\R;(L^\infty(\Rn\setminus0))^{m^2})$ and 
$\partial_t \mathcal{N}(t,\xi)\in 
L^1(\R;(L^\infty(\Rn\backslash0))^{m^2})$.
\end{lemma}

We have asymptotic integrations for (\ref{EQ:Eq})--(\ref{EQ:initial}).
\begin{proposition}[Matsuyama and Ruzhansky (\cite{MR-MN})] 
\label{prop:Rep}
Let $A(t,\xi)$ be a symbol of regularly hyperbolic operator $A(t,D_x)$ 
satisfying (\ref{hyp1}),  
and $\mathcal{N}(t,\xi)$ the diagonaliser of 
$A(t,\xi/|\xi|)$. 
Then there exist vector-valued functions 
${\bf a}^j(t,\xi)$, $j=0,1,\ldots,m-1$, 
determined by the initial value problem 
\[
D_t{\bf a}^j(t,\xi)=C(t,\xi){\bf a}^j(t,\xi), \qquad 
\left({\bf a}^1(0,\xi),\cdots,{\bf a}^m(0,\xi)\right)=\mathcal{N}(0,\xi),
\]
with
\[
C(t,\xi)=\Phi(t,\xi)^{-1} (D_t \mathcal{N}(t,\xi))
\mathcal{N}(t,\xi)^{-1}\Phi(t,\xi) \in L^1(\R; (L^\infty(\Rn\backslash0))^{m^2}),
\]
such that 
the solution $U(t,x)$ of (\ref{EQ:Eq})--(\ref{EQ:initial}) is represented by 
\begin{equation}\label{EQ:Representation of U}
U(t,x)=\sum_{j=0}^{m-1} 
\mathcal{F}^{-1} \left[\mathcal{N}(t,\xi)^{-1}\Phi(t,\xi)
{\bf a}^j(t,\xi)\widehat{f}_j(\xi) \right](x), 
\end{equation}
where we put 
\[
\Phi(t,\xi)
=\mathrm{diag}\left(e^{i\int_0^t\varphi_1(s,\xi)\,ds},\cdots,
e^{i\int_0^t\varphi_m(s,\xi)\,ds}\right).
\]
\end{proposition}

Let us introduce a class of symbols of differential operators, 
which is convenient for the fixed point argument. \\

\noindent 
{\bf Class $\mathcal{K}$.} {\em 
Given two constants $\Lambda>0$ and $K>0$, 
we say that a symbol $A(t,\xi)$    
of a pseudo-differential operator $A(t,D_x)$ 
belongs to 
$\mathcal{K}=\mathcal{K}(\Lambda,K)$ if $A(t,\xi/|\xi|)$ belongs to 
$\mathrm{Lip}_{\mathrm{loc}}(\R;(L^\infty(\Rn\backslash0))^{m^2})$ and satisfies
\[
\|A(t,\xi/|\xi|)\|_{L^\infty(\R;(L^\infty(\Rn\backslash0))^{m^2})} \le \Lambda,
\]
\[
\int^\infty_{-\infty}
\left\| \partial_t A(t,\xi/|\xi|)\right\|_{(L^\infty(\Rn\backslash0))^{m^2}}\, dt
\le K. 
\] 
}

The next lemma is the heart of our argument. It will be applied with a 
sufficiently small constant $K_0>0$ which will be fixed later.
\begin{lemma}\label{lem:Kirchhoff}
Let $n\ge1$. Assume that the symbol $A(t,\xi)$ 
of a pseudo-differential operator $A(t,D_x)$ satisfies strictly hyperbolic condition and 
(\ref{hyp1}),
and 
belongs to $\mathcal{K}=\mathcal{K}(\Lambda,K)$ 
for some $\Lambda>0$ and $0<K\le K_0$ with a sufficiently small constant 
$K_0>0$.
Let $U\in C(\R;{\bf L}^2(\Rn))$ 
be a solution to the Cauchy problem 
\[
D_tU=A(t,D_x)U, \quad U(0,x)=U_0(x)\in 
{\bf L}^2(\Rn) \cap \mathcal{Y}(\Rn),
\] 
and let $s(t)$ be the function
\[
s(t)=\langle SU(t,\cdot),U(t,\cdot)\rangle_{{\bf L}^2(\Rn)}. 
\]
Then there exist two constants $M>0$ and $c>0$
independent of $U$ and $K$  such that 
\begin{eqnarray}\label{core1}
& & \|A(s(t),\omega)\|_{(L^\infty({\bf S}^{n-1}))^{m^2}} \\
&\le& \|A(s(0),\omega)\|_{(L^\infty({\bf S}^{n-1}))^{m^2}}
+M \left(K\|U_0\|_{{\bf L}^2(\Rn)}^2
+\frac{1}{1-cK}\|U_0\|_{\mathcal{Y}(\Rn)}\right), 
\nonumber
\end{eqnarray}
\begin{eqnarray}\label{core2}
& & \int^\infty_{-\infty}\left\|\partial_t \left[A(s(t),\omega)\right] 
\right\|_{(L^\infty({\bf S}^{n-1}))^{m^2}} 
\, dt \\
&\le& M \left(K\|U_0\|_{{\bf L}^2(\Rn)}^2
+\frac{1}{1-cK}\|U_0\|_{\mathcal{Y}(\Rn)}\right).
\nonumber
\end{eqnarray} 
\end{lemma} 
{\em Outline of proof.} 
The proof is based on Lemma \ref{prop:Rep} via Fourier transform.
More prescisely, writing 
\[
s(t)=\langle S\widehat{U}(t,\xi),\widehat{U}(t,\xi)\rangle_{{\bf L}^2(\Rn)},
\]
we calculate the derivative $s^\prime(t)$ and plugging (\ref{EQ:Representation of U}) 
into $s^\prime(t)$. Then we can write 
\begin{equation}\label{EQ:s-prime}
s^\prime(t)=2 \mathrm{Re} \left\langle S\widehat{U}^\prime(t,\xi),
\widehat{U}(t,\xi) \right\rangle_{{\bf L}^2(\Rn)} 
=2\{I(t)+J(t)\}, 
\end{equation}
where
\begin{eqnarray*}
&& I(t)=\\
&& \mathrm{Re} \sum_{j,k=0}^{m-1}
\left\langle
S \mathcal{N}(t,\xi)^{-1}\partial_t \Phi(t,\xi){\bf a}^j(t,\xi)\widehat{f}_j(\xi),
\mathcal{N}(t,\xi)^{-1}\Phi(t,\xi){\bf a}^k(t,\xi)\widehat{f}_k(\xi)
\right\rangle_{{\bf L}^2(\Rn)},
\end{eqnarray*}
\begin{eqnarray*}
&& J(t)=\\
&& \mathrm{Re} \sum_{j,k=0}^{m-1}
\left\langle
S \pa_t \mathcal{N}(t,\xi)^{-1} \Phi(t,\xi){\bf a}^j(t,\xi)\widehat{f}_j(\xi),
\mathcal{N}(t,\xi)^{-1}\Phi(t,\xi){\bf a}^k(t,\xi)\widehat{f}_k(\xi)
\right\rangle_{{\bf L}^2(\Rn)}\\
&& +\left\langle
S \mathcal{N}(t,\xi)^{-1} \Phi(t,\xi)\partial_t {\bf a}^j(t,\xi)\widehat{f}_j(\xi),
\mathcal{N}(t,\xi)^{-1}\Phi(t,\xi){\bf a}^k(t,\xi)\widehat{f}_k(\xi)
\right\rangle_{{\bf L}^2(\Rn)}.
\end{eqnarray*}
We can proceed our analysis and conclude the proof.\\

\noindent 
{\em Outline of proof of Theorem \ref{thm:S-Kirchhoff}.}
We employ the Schauder-Tychonoff fixed point theorem. Let 
$A(t,\xi) \in \mathcal{K}$, and we fix 
the data $U_0\in {\bf L}^2(\Rn)\cap \mathcal{Y}(\Rn)$. Then it follows 
from Lemma~\ref{lem:Kirchhoff} that the mapping 
\[ \Theta : A(t,\xi)\mapsto A(s(t),\xi) \] 
maps $\mathcal{K}=\mathcal{K}(\Lambda,K)$ into itself provided that 
$\|U_0\|^2_{{\bf L}^2(\Rn)}+\|U_0\|_{\mathcal{Y}(\Rn)}$ 
is sufficiently small, with
$\Lambda>2\| A(0,\xi/|\xi|)\|_{(L^\infty(\Rn\backslash 0))^{m^2}}$
and sufficiently small $0<K<K_0$.
Now $\mathcal{K}$ may be regarded as the convex subset 
of the Fr\'echet space $L^{\infty}_{\mathrm{loc}}(\R;(L^\infty(\Rn\backslash0))^{m^2})$, 
and we endow 
$\mathcal{K}$ with the induced topology. 
We can show that $\mathcal{K}$ is 
compact in $L^{\infty}_{\mathrm{loc}}(\R;(L^\infty(\Rn\backslash0))^{m^2})$ 
and the mapping $\Theta$ is continuous on $\mathcal{K}$. 
Thus the Schauder-Tychonoff fixed point theorem allows us to conclude the proof.

%
%
%
%

\end{document}